\documentclass[11pt]{article}
\usepackage{latexsym}
\usepackage[all]{xy}
\usepackage{amsmath}
\usepackage{amssymb}
\usepackage{amsfonts}

\def\abstractname{R\'esum\'e}

\newtheorem{thm}{Th\'eor\`eme}

\newtheorem{lem}[thm]{Lemme}
\newtheorem{sublem}[thm]{Sous-lemme}
\newtheorem{subsublem}[thm]{Sous-sous-lemme}

\newtheorem{cor}[thm]{Corollaire}

\begin{document}

\title{\textbf{Anneaux de d\'efinition des dg-alg\`ebres
propres et lisses}}
\bigskip
\bigskip

\author{\bigskip\\
Bertrand To\"en\\
\small{Laboratoire Emile Picard} \\
\small{Universit\'e Paul Sabatier Toulouse 3, Bat 1R2}\\
\small{31 062 TOULOUSE cedex 9, France }\\
\small{e-mai:toen@math.ups-tlse.fr}}

\date{}

\maketitle

\begin{abstract}
Soit $k=colim_{i} k_{i}$ une colimite filtrante d'anneaux
commutatifs. On montre que la th\'eorie homotopique des
dg-alg\`ebres propres et lisses sur $k$ est la colimite 
des th\'eories homotopiques des
dg-alg\`ebres propres et lisses sur les $k_{i}$. Nous en d\'eduisons en particulier
que toute
dg-alg\`ebre propre et lisse est d\'efinissable sur une $\mathbb{Z}$-alg\`ebre
commutative de type fini. 
\end{abstract}

\medskip

\def\abstractname{Abstract}

\begin{abstract}
Let $k=colim_{i}k_{i}$ be a filtered colimit of commutative rings.
We show that the homotopy theory of smooth and proper dg-algebras
over $k$ is the colimit of the homotopy theories of smooth and proper
dg-algebras over the $k_{i}$'s. We deduce in particular that every 
smooth and proper dg-algebra can be defined over 
a commutative $\mathbb{Z}$-algebra of finite type.
\end{abstract}

\bigskip

\textbf{Classification MSC:} 16E45

\bigskip 

Ce texte accompagne le travail \cite{to}, et a pour objectif de d\'emontrer une autre 
propri\'et\'e de finitude sp\'ecifique aux dg-alg\`ebres propres et lisses. Cette propri\'et\'e
de finitude affirme qu'une telle dg-alg\`ebre sur une colimite filtrante 
d'anneaux commutatifs $A=colim A_{i}$ est toujours d\'efinissable sur un des $A_{i}$.
En particulier, une dg-alg\`ebre propre et lisse sur un anneau $k$
est toujours d\'efinissable sur une $\mathbb{Z}$-alg\`ebre commutative de type fini.
Ceci r\'epond   
positivement \`a la conjecture \cite[5.3]{ka}, et permet d'appliquer des techniques
de passage \`a la caract\'eristique positive pour l'\'etude des dg-alg\`ebres propres et lisses (dont
une application est l'\'etude de la d\'eg\'en\'erescence de la suite spectrale
de Hodge vers de Rham, voir \cite{ka}).

Les r\'esultats de ce travail sont dans la m\^eme lign\'ee que ceux de \cite{to}, bien qu'ils
ne semblent ni d\'ecouler de ceux-ci ni les impliquer directement. Tout comme nous l'avons mention\'e dans
\cite{to}, ces r\'esultats s'\'eclairent \`a la lumi\`ere de la th\'eorie des champs (sup\'erieurs et m\^eme 
d\'eriv\'es, voir \cite{to2}). Le th\'eor\`eme \ref{t1} d\'emontr\'e dans ce travail affirme en effet que le
champs des dg-alg\`ebres propres et lisses est localement de pr\'esentation finie. Ceci est 
un premier pas vers une preuve du fait que ce champ est un \emph{$D^{-}$-champ 
localement g\'eom\'etrique localement de pr\'esentation finie} (dans la terminologie de \cite{tova}). 
J'esp\`ere pouvoir revenir sur le caract\'ere g\'eom\'etrique du champ 
des dg-alg\`ebres propres et lisses dans un travail ult\'erieur, afin de mieux comprendre, et sans
doute d'am\'eliorer,  les r\'esultats de \cite{to}. \\

\textbf{Remerciements:} Je tiens \`a remercier D. Kaledin pour de nombreuses
conversations sur le sujet, et en particulier pour avoir remarquer 
que le sous-lemme \ref{slkar} d\'ecoulait formellement des arguments 
de la preuve du th\'eor\`eme principal. Je remercie aussi le rapporteur pour
ses commentaires. \\

\bigskip

Nous renvoyons \`a \cite{to,tova,ke,koso} pour les notions de dg-alg\`ebres propres, lisses et homotopiquement 
de pr\'esentation finie. Pour un anneau commutatif $k$ nous noterons 
$Ho(k-dg-alg^{pl})$ la sous-cat\'egorie pleine de la cat\'egorie homotopique des $k$-dg-alg\`ebres
form\'ee des dg-alg\`ebres propres et lisses sur $k$. Pour tout morphisme d'anneaux
commutatifs $k \longrightarrow k'$, on dispose d'un foncteur de changement de
bases 
$$-\otimes_{k}^{\mathbb{L}}k' : Ho(k-dg-alg) \longrightarrow Ho(k'-dg-alg).$$
Ce foncteur pr\'eserve les caract\`eres propres, lisses, et homotopiquement
de pr\'esentation finie des dg-alg\`ebres. Avec ces notations, le th\'eor\`eme principal 
de ce travail est le suivant.

\begin{thm}\label{t1}
Soit $I$ une cat\'egorie filtrante et $\{k_{i}\}_{i\in I}$ un $I$-diagramme
d'anneaux commutatif de colimite $k:=colim_{i\in I}k_{i}$. Le foncteur naturel
$$Colim_{i\in I}-\otimes^{\mathbb{L}}_{k_{i}}k : 
Colim_{i\in I} Ho(k_{i}-dg-alg^{pl}) \longrightarrow Ho(k-dg-alg^{pl})$$
est une \'equivalence de cat\'egories. 
\end{thm}

La preuve que nous donnons de ce th\'eor\`eme repose sur les faits connus suivants, tous d\'emontr\'es
dans \cite{tova} \`a l'exception de F7 dont nous donnerons une preuve (voir le sous-lemme \ref{slkar}).

\begin{enumerate}

\item[F1]\label{f1} Une dg-alg\`ebre propre et lisse est 
homotopiquement de pr\'esentation finie (voir \cite[Cor. 2.13]{tova}).

\item[F2]\label{f2} Une dg-alg\`ebre homotopiquement de pr\'esentation finie est 
lisse (voir \cite[Prop. 2.14]{tova}).

\item[F3]\label{f3} Une dg-alg\`ebre est homotopiquement de pr\'esentation finie si et seulement si
elle est quasi-isomorphe \`a un r\'etracte d'un complexe $I$-cellulaire fini (voir \cite[Prop. 2.2]{tova}). 

\item[F4]\label{f4} Soit $I$ une cat\'egorie filtrante et $\{k_{i}\}_{i\in I}$ un $I$-diagramme
d'anneaux commutatifs de colimite $k:=colim_{i\in I}k_{i}$. Soit $A_{i}$ une 
$k_{i}$-dg-alg\`ebre (pour un $i\in I$) et notons
$$A:=A_{i}\otimes_{k_{i}}^{\mathbb{L}}k\simeq Colim_{j\in i/I}A_{i}\otimes_{k_{i}}^{\mathbb{L}}k_{j}.$$
Alors le foncteur naturel 
$$Colim_{j\in i/I}-\otimes^{\mathbb{L}}_{k_{j}}k : 
Colim_{j\in i/I}D_{parf}(A_{i}\otimes_{k_{i}}^{\mathbb{L}}k_{j}) \longrightarrow
D_{parf}(A)$$
est une \'equivalence de cat\'egories (voir \cite[Lem. 2.10]{tova}). 

\item[F5]\label{f5} Soit $I$ une cat\'egorie filtrante et $\{k_{i}\}_{i\in I}$ un $I$-diagramme
d'anneaux commutatifs de colimite $k:=colim_{i\in I}k_{i}$. Soit $A_{i}$ une 
$k_{i}$-dg-alg\`ebre (pour un $i\in I$) homotopiquement de pr\'esentation finie
et notons
$$A:=A_{i}\otimes_{k_{i}}^{\mathbb{L}}k\simeq Colim_{j\in i/I}A_{i}\otimes_{k_{i}}^{\mathbb{L}}k_{j}.$$
Nous noterons $D_{pspa}(A)$ (resp. $D(A_{i}\otimes_{k_{i}}^{\mathbb{L}}k_{j})$)
la sous-cat\'egorie pleine de $D(A)$ (resp. de $D_{pspa}(A_{i}\otimes_{k_{i}}^{\mathbb{L}}k_{j})$)
form\'ee des $A$-dg-modules (resp. des $A_{i}\otimes_{k_{i}}^{\mathbb{L}}k_{j}$-dg-modules) 
dont le complexe sous-jacent de $k$-modules (resp. de $k_{j}$-modules) est parfait. 
Alors le foncteur naturel 
$$Colim_{j\in i/I}-\otimes^{\mathbb{L}}_{k_{j}}k : 
Colim_{j\in i/I}D_{pspa}(A_{i}\otimes_{k_{i}}^{\mathbb{L}}k_{j}) \longrightarrow
D_{pspa}(A)$$
est une \'equivalence de cat\'egories. Ceci est une cons\'equence du caract\`ere
localement de pr\'esentation finie du champ $\mathcal{M}_{T}$ de \cite[Thm. 3.6]{tova}
(pour $T=A_{i}$). Par soucis de clart\'e nous redonnerons une preuve de
ce fait (voir le sous-lemme \ref{sl}).

\item[F6]\label{f6} Soit $A$ une $k$-dg-alg\`ebre lisse. Alors un $A$-dg-module
$E$ qui est parfait comme complexe de $k$-modules et aussi parfait
comme $A$-dg-module (voir \cite[Lem. 2.8 (2)]{tova}). En d'autres
termes on a $D_{pspa}(A)\subset D_{parf}(A)$. 

\item[F7]\label{f7} La cat\'egorie homotopique $Ho(k-dg-alg)$ est 
Karoubienne. En d'autres termes, pour tout objet $A\in Ho(k-dg-alg)$
et tout endomorphisme idempotent $p$ de $A$, il existe 
un objet $B$ et deux morphismes $i : B \longrightarrow A$, $r : A \longrightarrow B$ tels
que $ri=id$ et $ir=p$. Ce fait n'est pas une formalit\'e car il est bien connu 
qu'il existe des cat\'egories de mod\`eles $M$ telle que $Ho(M)$ ne soit 
pas Karoubienne (l'exemple le plus fameux est celui de la cat\'egorie homotopique
des espaces, voir par exemple \cite{fh}).

\end{enumerate}

Passons maintenant \`a la preuve du th\'eor\`eme \ref{t1}. 
Nous commencerons par d\'emontrer l'analogue du th\'eor\`eme \ref{t1} pour 
les dg-alg\`ebres homotopiquement de pr\'esentation finie. Pour cela, nous noterons
$Ho(k-dg-alg^{hpf})$ la sous-cat\'egorie pleine de $Ho(k-dg-alg)$ form\'ee
des $k$-dg-alg\`ebres homotopiquement de pr\'esentation finie (pour un
anneau $k$). 

\begin{lem}\label{l1}
Soit $I$ une cat\'egorie filtrante et $\{k_{i}\}_{i\in I}$ un $I$-diagramme
d'anneaux commutatifs de colimite $k:=colim_{i\in I}k_{i}$. Le foncteur naturel
$$Colim_{i\in I}-\otimes^{\mathbb{L}}_{k_{i}}k : 
Colim_{i\in I} Ho(k_{i}-dg-alg^{hpf}) \longrightarrow Ho(k-dg-alg^{hpf})$$
est une \'equivalence de cat\'egories. 
\end{lem}

\textit{Preuve du lemme \ref{l1}:} Soient 
$A_{i}$ et $B_{i}$ deux objets de $Ho(k_{i}-dg-alg^{hpf})$ pour un objet $i\in I$. 
Notons 
$$A:=A\otimes_{k_{i}}^{\mathbb{L}}k \qquad B:=B\otimes_{k_{i}}^{\mathbb{L}}k.$$
Vue comme $k_{i}$-dg-alg\`ebre (\`a travers le foncteur d'oubli), la dg-alg\`ebre
$B$ s'\'ecrit comme une colimite filtrante 
$$B\simeq Colim_{j\in i/I}B\otimes_{k_{i}}^{\mathbb{L}}k_{j}.$$
Ainsi, comme $A$ est homotopiquement de pr\'esentation finie on trouve 
$$[A,B]_{Ho(k-dg-alg)}\simeq
[A_{i},Colim_{j\in i/I}B\otimes_{k_{i}}^{\mathbb{L}}k_{j}]_{Ho(k_{i}-dg-alg)}$$
$$\simeq 
Colim_{j\in i/I}[A_{i},B\otimes_{k_{i}}^{\mathbb{L}}k_{j}]_{Ho(k_{i}-dg-alg)}
\simeq Colim_{j\in i/I}[A\otimes_{k_{i}}^{\mathbb{L}}k_{j},B\otimes_{k_{i}}^{\mathbb{L}}k_{j}]_{Ho(k_{j}-dg-alg)}.$$
Ceci montre que le foncteur en question est pleinement fid\`ele.

Pour l'essentielle surjectivit\'e nous allons utiliser F3 et F7. Soit $A\in Ho(k-dg-alg^{hpf})$. D'apr\`es
F3 nous savons qu'il existe une $k$-dg-alg\`ebre
$B$ qui est quasi-isomorphe \`a un complexe $I$-cellulaire fini, et deux morphismes
$$\xymatrix{A \ar[r]^-{i} & B \ar[r]^-{r} & A,}$$
tels que $r\circ i=id$ (dans $Ho(k-dg-alg^{hpf})$). Notons $p=i\circ r$, qui est un 
endomorphisme de $B$ tel que $p^{2}=p$ (toujours dans $Ho(k-dg-alg^{hpf})$).  

Comme $B$ est un complexe $I$-cellulaire fini, il est facile de voir que $B$ est d\'efinie, \`a quasi-isomorphisme
pr\`es, 
sur $k_{i}$ pour $i\in I$. Soit donc $B_{i}$, une $k_{i}$-dg-alg\`ebre, qui est aussi quasi-isomorphe \`a un
complexe $I$-cellulaire fini, telle que $B_{i}\otimes_{k_{i}}^{\mathbb{L}}k\simeq B$ dans $Ho(k-dg-alg)$. 
On peut donc \'ecrire
$$B\simeq Colim_{j\in i/I}B_{i}\otimes_{k_{i}}^{\mathbb{L}}k_{j}.$$
Comme $B$ est un
complexe $I$-cellulaire fini c'est un objet homotopiquement de pr\'esentation finie. 
On a donc,
$$[B,B]_{Ho(k-dg-alg)}\simeq Colim_{j\in i/I}
[B_{i}\otimes_{k_{i}}k_{j},B_{i}\otimes_{k_{i}}k_{j}]_{Ho(k_{j}-dg-alg)}.$$
Ainsi, quitte \`a remplacer l'objet $i$ par un objet $j\in i/I$, on peut supposer qu'il 
existe $p_{i}$ un endomorphisme de $B_{i}$ dans $Ho(k_{i}-dg-alg)$, tel que
$$p_{i}^{2}=p_{i} \qquad p_{i}\otimes_{k_{i}}^{\mathbb{L}}k=p.$$

Nous invoquons ici le fait F7, qui d\'ecoule du sous-lemme plus g\'en\'eral suivant
(en prenant $M=k-dg-alg$, et $\Pi:=H^{*}$ le foncteur de cohomologie totale).

\begin{sublem}\label{slkar}
Soit $M$ une cat\'egorie de mod\`eles munie d'un foncteur 
$\Pi : M \longrightarrow Ens$ v\'erifiant les propri\'et\'es suivantes.
\begin{enumerate}
\item Un morphisme $f : A \longrightarrow B$ est une \'equivalence si et seulement si
le morphisme induit $\Pi (f) : \Pi (A) \longrightarrow \Pi (B)$
est un isomorphisme.
\item Le foncteur $\Pi$ commute aux colimites filtrantes.
\end{enumerate}
Alors la cat\'egorie $Ho(M)$ est Karoubienne.
\end{sublem}

Avant de donner la preuve du sous-lemme pr\'ec\'edent montrons comment il permet 
de conclure la preuve du lemme \ref{l1}. On choisit un objet $A_{i}$ dans
$Ho(k_{i}-dg-alg)$ ainsi que deux morphismes
$i_{i} : A_{i} \longrightarrow B_{i}$ et $r_{i} : B_{i} \longrightarrow A_{i}$ avec 
$i_{i}r_{i}=p_{i}$ et $r_{i}i_{i}=id$. Alors $A_{i}\otimes_{k_{i}}^{\mathbb{L}}k$ et $A$ sont
isomorphes dans $Ho(k-dg-alg)$. En effet, cela se d\'eduit de l'unicit\'e de la
d\'ecomposition d'un endomorphisme idempotent. Plus pr\'ecis\`emment, on consid\`ere 
les deux morphismes
$$\xymatrix{
A_{i}\otimes_{k_{i}}^{\mathbb{L}}k \ar[r]^-{i_{i}\otimes_{k_{i}}^{\mathbb{L}}k} & B_{i}\otimes_{k_{i}}^{\mathbb{L}}k\simeq B \ar[r]^-{r} & 
A}$$
$$\xymatrix{
A \ar[r]^-{i} & B\simeq B_{i}\otimes_{k_{i}}^{\mathbb{L}}k \ar[r]^-{r_{i}\otimes_{k_{i}}^{\mathbb{L}}k} & A_{i}.}$$
Il est facile de v\'erifier que ces deux morphismes sont inverses l'un de l'autre. \\

\textit{Preuve du sous-lemme \ref{slkar}:} 
Soit $B$ un objet de $M$ et $p$ un endomorphisme idempotent de
$B$ dans $Ho(M)$. On peut supposer que $B$ est un objet cofibrant et fibrant et donc
repr\'esenter $p$ par un morphisme $e : B \longrightarrow B$
dans $M$. Nous d\'efinissons alors un objet
$$A:=Colim(\xymatrix{
B \ar[r]^-{e} & B \ar[r]^-{e} & \dots,})$$
qui est la composition transfinie d\'enombrable du morphisme $e$. 

Nous allons maintenant construire deux morphismes naturels
$$i : A \longrightarrow B \qquad r : B \longrightarrow A,$$
et v\'erifier que $ir=p$ et $ri=id$.
Pour cela, notons $\mathbb{N}$ la cat\'egorie associ\'ee \`a l'ensemble
ordon\'e des entiers naturels. Nous disposons de la cat\'egorie des $\mathbb{N}$-diagrammes
$M^{\mathbb{N}}$. Cette cat\'egorie est elle m\^eme une 
cat\'egorie de mod\`ele pour la structure de Reedy (voir par exemple \cite{ho}), 
pour laquelle les fibrations et \'equivalences sont d\'efinies termes \`a termes. 
Nous diposons d'un foncteur naturel
$$Ho(M^{\mathbb{N}}) \longrightarrow Ho(M)^{\mathbb{N}}$$
de la cat\'egorie homotopique des $\mathbb{N}$-diagrammes dans 
$M$ vers celle des $\mathbb{N}$-diagrammes dans $Ho(M)$. 

\begin{subsublem}\label{sl0}
Le foncteur ci-dessus 
$$Ho(M^{\mathbb{N}}) \longrightarrow Ho(M)^{\mathbb{N}}$$
est plein.
\end{subsublem}

\textit{Preuve du sous-sous-lemme \ref{sl0}:} Soient 
$$X:=(\xymatrix{X_{0} \ar[r]^-{f_{0}} & X_{1} \ar[r]^-{f_{1}} & \dots & X_{n} \ar[r]^-{f_{n}} & X_{n+1} \ar[r] & \dots})$$
$$Y:=(\xymatrix{Y_{0} \ar[r]^-{g_{0}} & Y_{1} \ar[r]^-{g_{1}} & \dots & Y_{n} \ar[r]^-{g_{n}} & Y_{n+1} \ar[r] & \dots})$$
deux $\mathbb{N}$-diagrammes dans $Ho(M^{\mathbb{N}})$. On supposera 
$X$ cofibrant et $Y$ fibrant. Cela signifie que chaque $Y_{n}$ est fibrant 
dans $M$, et que chaque morphisme $f_{n} : X_{n} \rightarrow X_{n+1}$ est 
une cofibration entre objets cofibrants.

Soit $u : X \rightarrow Y$ un morphisme dans $Ho(M)^{\mathbb{N}}$. Le morphisme
$u$ est donn\'e par des morphismes $u_{n} \in [X_{n},Y_{n}]$ dans $Ho(M)$, tels que
pour tout $n$ on ait $g_{n} u _{n}=u_{n+1}f_{n}$ dans $Ho(M)$. Comme
$X_{n}$ est cofibrant et $Y_{n}$ est fibrant, on peut choisir un 
repr\'esentant $u'_{n} : X_{n} \longrightarrow Y_{n}$ dans $M$. L'\'egalit\'e pr\'ec\'edente 
implique alors que les deux morphismes
$g_{n} u' _{n}$ et $\; u'_{n+1} f_{n}$ 
sont homotopes. Par r\'ecurrence sur $n$, nous allons construire
des morphismes $v_{n} : X_{n} \longrightarrow Y_{n}$, qui sont
homotopes aux $u'_{n}$, et tels que $g_{n} v _{n}=v_{n+1}f_{n}$
dans $M$. Supposons que de tels morphismes $v_{i}$ soient 
construits pour $i\leq n$. On choisit une homotopie
$$\xymatrix{ & Y_{n+1} & \\ 
 X_{n} \ar[rd]_-{u'_{n+1}f_{n}} \ar[ru]^-{g_{n} v _{n}}\ar[r]^-{h} & \ar[d]^-{p_{1}} \ar[u]_-{p_{0}}  Cyl(Y_{n+1}) \\
 & Y_{n+1}  }
$$
o\`u $Cyl(Y_{n})$ est un objet cylindre pour $Y_{n+1}$ (voir \cite{ho}). 
On consid\`ere le diagramme suivant
$$\xymatrix{
X_{n} \ar[d]_-{f_{n}} \ar[r]^-{h} & Cyl(Y_{n+1}) \ar[d]^-{p_{1}} \\
X_{n+1} \ar[r]_-{u'_{n+1}} & Y_{n+1}.}$$
Le morphisme $p_{1}$ \'etant une fibration triviale on peut choisir 
un rel\`evement $k : X_{n+1} \longrightarrow Cyl(Y_{n+1})$ de $h$ (qui 
fasse commuter les deux triangles possibles). On pose
$$v_{n+1}:=p_{0}k : X_{n+1} \longrightarrow Y_{n+1}.$$
Par construction, $v_{n+1}$ est homotope \`a $u'_{n+1}$ (car 
$k$ est une homotopie), et $v_{n+1}f_{n}=g_{n}v_{n}$. 
Les morphismes $v_{n}$ d\'efinissent alors un morphisme
$X \longrightarrow Y$ dans $M^{\mathbb{N}}$, dont l'image
dans $Ho(M)^{\mathbb{N}}$ est le morphisme
$u$ que l'on s'est donn\'e. 
\hfill $\Box$ \\

Revenons maintenant \`a la preuve du sous-lemme \ref{slkar}.
On consid\`ere
deux objets 
$$X:=(\xymatrix{B \ar[r]^-{e} & B \ar[r]^-{e} &  \dots B \ar[r]^-{e} & B \ar[r] & \dots})$$
$$Y:=(\xymatrix{B \ar[r]^-{id} & B \ar[r]^-{id} & \dots  B \ar[r]^-{id} & B \ar[r] & \dots})$$
dans $Ho(M^{\mathbb{N}})$. On dispose d'un morphisme 
$X \longrightarrow Y$ dans $Ho(M)^{\mathbb{N}}$ d\'efini par 
$u_{n}:=e$ pour tout $n$
$$\xymatrix{
B \ar[r]^-{e} \ar[d]^-{e} & \ar[d]^-{e} B \ar[r]^-{e} & 
\dots  \ar[d]^-{e} B \ar[r]^-{e} & \ar[d]^-{e} B \ar[r] & \dots \\
B \ar[r]^-{id} & B \ar[r]^-{id} & \dots   B \ar[r]^-{id} & B \ar[r] & \dots.}$$
Par le sous-sous-lemme \ref{sl0} ce morphisme se rel\`eve en un morphisme
$i' : X \longrightarrow Y$ dans $Ho(M^{\mathbb{N}})$.  De m\^eme, on dispose d'un morphisme 
$Y \longrightarrow X$ dans $Ho(M)^{\mathbb{N}}$ en posant $u_{n}:=e$
pour tout $n$. Ce morphisme se rel\`eve en un morphisme $r' : Y \longrightarrow X$
dans $Ho(M^{\mathbb{N}})$.
Finalement, en appliquant le foncteur (notons que les propri\'et\'es
$(1)$ et $(2)$ du sous-lemme \ref{slkar} impliquent que les colimites filtrantes
pr\'es\`ervent les \'equivalences dans $M$)
$$Colim_{\mathbb{N}} : Ho(M)^{\mathbb{N}} \longrightarrow Ho(M)$$
on trouve nos deux morphismes $i$ et $r$ dans $Ho(M)$
$$i:=colim(i') : A \simeq colim(X) \longrightarrow B \simeq colim(Y)$$
$$r:=colim(r') : B \simeq colim(Y) \longrightarrow A \simeq colim(X).$$

Consid\'erons maintenant $u : r\circ i : A \longrightarrow A$ dans $Ho(M)$. 
Cet endomorphisme induit un morphisme 
$$\Pi (u) : \Pi (A) \longrightarrow \Pi (A).$$
Comme $\Pi$  commute aux colimites filtrantes 
$\Pi (u)$ est \'egal au morphisme compos\'e
$$Colim (\xymatrix{
\Pi (B) \ar[r]^-{\Pi (e)} & \Pi (B) \ar[r]^-{\Pi (e)} & \dots,}) \longrightarrow 
Colim (\xymatrix{
\Pi (B) \ar[r]^-{id} & \Pi (B) \ar[r]^-{id} & \dots,})$$
$$ \longrightarrow Colim (\xymatrix{
\Pi (B) \ar[r]^-{\Pi (e)} & \Pi (B) \ar[r]^-{\Pi (e)} & \dots,})$$
tous deux induits par $\Pi (e)$ composantes par composantes.
Enfin, comme $\Pi (e)$ est un endomorphisme idempotent de $\Pi (B)$ il n'est pas
difficile de voir que $\Pi (u)=id$. Ceci implique en particulier que 
$u$ est un automorphisme de $A$ dans $Ho(M)$.

Pour conclure, on remarque que par construction nous avons
$i\circ r=p$ (dans $Ho(M)$), comme on peut le voir en consid\'erant le diagramme
commutatif suivant
$$\xymatrix{
B \ar[r]^-{\sim} \ar[d]_-{e} & Colim (Y) \ar[d]^-{r} \\
B \ar[r] \ar[d]_-{e} & Colim(X) \ar[d]^-{i} \\
B \ar[r]^-{\sim} & Colim(Y). }$$
Ainsi, on a 
$$u^{3}=r\circ p^{2}\circ i=r\circ p \circ i=u^{2}.$$
Comme nous avons vu que $u=r\circ i$ est un automorphisme cela implique
que $r\circ i=id$. 
 \hfill $\Box$ \\

Ceci termine la preuve du lemme \ref{l1}.
\hfill $\Box$ \\

Un corollaire imm\'ediat du lemme \ref{l1} est le r\'esultat suivant.

\begin{cor}\label{c1}
Soit $I$ une cat\'egorie filtrante et $\{k_{i}\}_{i\in I}$ un $I$-diagramme
d'anneaux commutatifs de colimite $k:=colim_{i\in I}k_{i}$. Le foncteur naturel
$$Colim_{i\in I}-\otimes^{\mathbb{L}}_{k_{i}}k : 
Colim_{i\in I} Ho(k_{i}-dg-alg^{pl}) \longrightarrow Ho(k-dg-alg^{pl})$$
est pleinement fid\`ele
\end{cor}

\textit{Preuve:} Cela se d\'eduit formellement de F1 et du lemme \ref{l1}. \hfill $\Box$

\begin{lem}\label{l2}
Soit $I$ une cat\'egorie filtrante et $\{k_{i}\}_{i\in I}$ un $I$-diagramme
d'anneaux commutatifs de colimite $k:=colim_{i\in I}k_{i}$. Soit 
$A_{i}$ une $k_{i}$-dg-alg\`ebre homotopiquement de pr\'esentation finie pour un $i\in I$. 
Si $A_{i}\otimes_{k_{i}}^{\mathbb{L}}k$ est une $k$-dg-alg\`ebre propre, alors
il existe $j\in i/I$ tel que $A_{i}\otimes_{k_{i}}^{\mathbb{L}}k_{j}$ soit 
une $k_{j}$-dg-alg\`ebre propre.
\end{lem}

\textit{Preuve:} Soit $A:=A_{i}\otimes_{k_{i}}^{\mathbb{L}}k$. On consid\`ere
la $A$-dg-module parfait $A\in D_{parf}(A)$. Nous souhaitons alors utiliser le
fait F5 afin de montrer que ce dg-module est d\'efini sur $A_{j}$ pour un $j\in i/I$. 
Avant cela nous red\'emontrerons ce fait dans le sous-lemme suivant. 

\begin{sublem}\label{sl}
Soit $I$ une cat\'egorie filtrante et $\{k_{i}\}_{i\in I}$ un $I$-diagramme
d'anneaux commutatifs de colimite $k:=colim_{i\in I}k_{i}$. Soit $A_{i}$ une 
$k_{i}$-dg-alg\`ebre (pour un $i\in I$) homotopiquement de pr\'esentation finie
et notons
$$A:=A_{i}\otimes_{k_{i}}^{\mathbb{L}}k\simeq Colim_{j\in i/I}A_{i}\otimes_{k_{i}}^{\mathbb{L}}k_{j}.$$
Alors le foncteur naturel 
$$Colim_{j\in i/I}-\otimes^{\mathbb{L}}_{k_{j}}k : 
Colim_{j\in i/I}D_{pspa}(A_{i}\otimes_{k_{i}}^{\mathbb{L}}k_{j}) \longrightarrow
D_{pspa}(A)$$
est une \'equivalence de cat\'egories. 
\end{sublem}

\textit{Preuve du sous-lemme:} La pleine fid\`elit\'e du foncteur en question se d\'eduit ais\`emment 
des faits F2, F4 et F6.  Il nous faut montrer l'essentielle surjectivit\'e. Pour cela, soit
$E\in D_{pspa}(A)$. Nous pouvons bien entendu supposer que $A$ est une
dg-alg\`ebre cofibrante, et que $E$ est un $A$-dg-module cofibrant. 
Ceci implique en particulier que le complexe sous-jacent \`a $E$, que nous
noterons $M$, est un objet cofibrant dans la cat\'egorie de mod\`eles
des complexes de $k$-modules. 
Le $A$-dg-module $E$ est alors determin\'e (voir par exemple \cite[Lem. 3.14]{tova})
par la classe d'isomorphisme de l'objet $M\in D_{parf}(k)$, et par le morphisme 
$A \longrightarrow \underline{End}(M)$ dans $Ho(k-dg-alg)$, 
o\`u $\underline{End}(M)$ est la $k$-dg-alg\`ebre des endomorphismes du complexe $M$. 
Comme $M$ est un complexe parfait, la propri\'et\'e F4 (appliqu\'ee \`a $A=k$) implique
qu'il existe $j\in i/I$ et $M_{j} \in D_{parf}(k_{j})$ tel que 
$M_{j}\otimes_{k_{j}}^{\mathbb{L}}k\simeq M$. Nous pouvons supposer que $i=j$. 
Dans ce cas, on a un isomorphisme dans $Ho(k_{i}-dg-alg)$
$$\underline{End}(M)\simeq Colim_{j\in i/I}\underline{End}(M_{i}\otimes_{k_{i}}^{\mathbb{L}}k_{j})$$
(en supposant que $M_{i}$ a \'et\'e choisi cofibrant). Comme 
$A_{i}$ est homotopiquement de pr\'esentation finie, cela implique que le morphisme
$$A_{i} \longrightarrow \underline{End}(M)$$
se factorise dans $Ho(k_{i}-dg-alg)$ par un morphisme
$$A_{i} \longrightarrow \underline{End}(M_{j}).$$
Ce morphisme induit un morphisme dans $Ho(k_{j}-dg-alg)$
$$A_{i}\otimes_{k_{i}}^{\mathbb{L}}k_{j} \longrightarrow \underline{End}(M_{j}).$$
Ce dernier morphisme determine un $A_{i}\otimes_{k_{i}}^{\mathbb{L}}k_{j}$-dg-module
$E_{j} \in D(A_{i}\otimes_{k_{i}}^{\mathbb{L}}k_{j})$
qui par construction est tel que $E_{j}\otimes_{k_{j}}^{\mathbb{L}}k \simeq E$. 
Comme $E_{j}$ est parfait comme complexe de $k_{j}$-module, ceci montre l'essentielle
surjectivit\'e du foncteur en question. \hfill $\Box$ \\

Nous revenons maintenant \`a la preuve du lemme \ref{l2}.
Par hypoth\`ese de propret\'e sur $A$, le  $A$-dg-module $A$ est 
parfait comme complexe de $k$-modules, et d'apr\`es le sous-lemme \ref{sl}, quitte \`a 
remplacer $i$ par un $j\in i/I$, on peut supposer qu'il 
existe $E_{i}\in D_{pspa}(A_{i})$ tel que $E_{i}\otimes_{k_{i}}^{\mathbb{L}}k\simeq A$
dans $D(A)$. Or $A_{i}$ est homotopiquement de pr\'esentation finie et donc lisse d'apr\`es F2. 
En particulier, F6 implique que $E_{i}$ est un $A_{i}$-dg-module parfait. Ainsi, 
$E_{i}$ et $A_{i}$ sont deux $A_{i}$-dg-modules parfaits qui deviennent
quasi-isomorphes apr\`es extension des scalaires \`a $k$. La propri\'et\'e F4 implique
donc  qu'il existe $j\in i/I$ tel que $E_{i}\otimes_{k_{i}}^{\mathbb{L}}k_{j}\simeq
A_{i}\otimes_{k_{i}}^{\mathbb{L}}k_{j}$ dans $D(A_{i}\otimes_{k_{i}}^{\mathbb{L}}k_{j})$. 
Cela a comme cons\'equence que $A_{i}\otimes_{k_{i}}^{\mathbb{L}}k_{j}$ est un 
complexe parfait de $k_{j}$-modules, et donc est une $k_{j}$-dg-alg\`ebre
propre. \hfill $\Box$ \\

Nous pouvons maintenant terminer la preuve du th\'eor\`eme \ref{t1}. D'apr\`es \ref{c1} le foncteur
en question est pleinement fid\`ele. Soit alors $A$ une $k$-dg-alg\`ebre propre et lisse.
D'apr\`es \ref{l2}, il existe $i\in I$ et une $k_{i}$-dg-alg\`ebre homotopiquement de pr\'esentation finie $A_{i}$ avec 
$A_{i}\otimes_{k_{i}}^{\mathbb{L}}k\simeq A$. Le lemme \ref{l2} nous dit alors qu'il existe 
$j\in i/I$ tel que $A_{j}:=A_{i}\otimes_{k_{i}}^{\mathbb{L}}k_{j}$ soit une
$k_{j}$-dg-alg\`ebre propre et lisse avec 
$A_{j}\otimes_{k_{j}}^{\mathbb{L}}k$ quasi-isomorphe \`a $A$. 
Ceci montre que le foncteur en question dans le th\'eor\`eme \ref{t1} est essentiellement surjectif, et 
d\'emontre donc ce th\'eor\`eme en vu du corollaire \ref{c2}. \\

Une cons\'equence imm\'ediate du th\'eor\`eme \ref{t1} est la corollaire suivant, qui est une r\'eponse 
positive \`a la conjecture \cite[5.3]{ka}. 

\begin{cor}\label{c2}
Soit $k$ un anneau commutatif et $A$ une $k$-dg-alg\`ebre propre et lisse. Alors il existe
un anneau commutatif $k_{0}$, de type fini sur $\mathbb{Z}$, un morphisme 
$k_{0} \longrightarrow k$, et une $k_{0}$-dg-alg\`ebre propre et lisse $A_{0}$ telle que
$$A_{0}\otimes_{k_{0}}^{\mathbb{L}}k\simeq A.$$
\end{cor}

\textit{Preuve:} On \'ecrit $k$ comme une colimite filtrante de sous-anneaux de type
fini sur $\mathbb{Z}$, et on applique le th\'eor\`eme \ref{t1}. \hfill $\Box$ \\

\end{document}